\documentclass[twoside,leqno,11pt]{article}
\usepackage{ifthen,amsmath,graphicx,epsfig}
\numberwithin{equation}{section} \numberwithin{lemma}{section}
\numberwithin{theorem}{section} \numberwithin{corollary}{section}
\numberwithin{proposition}{section}
\numberwithin{definition}{section} \numberwithin{example}{section}
\numberwithin{remark}{section}
\def\TITLE       {ON THE DIOPHANTINE EQUATION $x^{4}-q^{4}=py^{5}$}
\medskip
\def\LASTNAMEI   {Savin}
\def\FIRSTNAMEI  {Diana}
\def\ABSTRACT {
\begin{center}
ABSTRACT:
\end{center}
In this paper we study the Diophantine equation $x^{4}-q^{4}=py^{5},$ 
with the following conditions: $p$ and $q$ are
different prime natural numbers, $y$ is not divisible with $p$, $p\equiv3$ (mod$20$), $q\equiv4$ (mod$5$), $\overline{p}$ is a generator of the group $\left(U\left(\textbf{Z}_{q^{4}}\right),\cdot\right)$, $\left(x,y\right)=1$, $2$ is a $5$-power residue mod $q$.\ \ \

}
\def\CLASSIFICATION{MSC (2000): 11D41}

\def\KEYWORDS{KEYWORDS: Diophantine equations}

\begin{document}

\begin{center}
\par\TITLE
\par\FIRSTNAMEI \ \LASTNAMEI
\par\medskip

\end{center}
\par\
\par\ABSTRACT{}
\par\
\par\CLASSIFICATION{}
\par\KEYWORDS{}
\par\


 \section{Introduction}

In some previous papers, [5], [6], [8] we have solved Diophantine equations of the form\\
$$x^{4}-y^{4}=pz^{2},$$\\
where $p$ is a prime natural number from the set $\left\{3,5,7,11,13,19,29,37\right\}$.\\
In the paper[9], we have solved Diophantine equation of the form\\
$$x^{4}-q^{4}=py^{3},$$\\
where $p$ and $q$ are prime natural numbers, with the special conditions.\\
In this paper we continue our study on special Diophantine equations of order $4$, by considering an equation of order $5$, with special coeffiecients.\\
The study is done by using results on rings of algebraic integers associated in Kummer fields.\\ 
A special use is considered for the Hilbert's symbol $\left\{\frac{\alpha}{P}\right\}$ for the number $\alpha$ relatively to the prime ideal $P$ in the ring $Z\left[\xi\right]$, given in his book "Theory of algebraic number fields", which we have read in Romanian translation.\\
It is interesting how old results in theory of numbers could lead us to new results in theory of Diophantine equations.\\
\bigskip\\
In the begining we recall some results.  
\bigskip\\
\textbf{Proposition 1.1. }( $\left[2\right] $).
\textit{Let} $l$ \textit{be a natural number} $l\geq3$ \textit{and} $\xi$ \textit{be a primitive root of unity of order} $l$, $Z\left[\xi\right]$ \textit{be the ring of integers of the cyclotomic field} $\textbf{Q}\left(\xi\right)$. \textit{If} $p$ \textit{is a prime natural number}, $l$ \textit{is not divisible with} $p$, \textit{and} $f$ \textit{is the smallest positive integer such that} $p^{f}$$\equiv$$1$ (\textit{mod}$l$), \textit{then we have}:\\
$$pZ\left[\xi\right]=P_{1}P_{2}...P_{r},$$\\
 \textit{where} $r=\frac{\varphi\left(l\right)}{f}$, $P_{j}$ $j=\overline{1,r}$
 \textit{are different prime ideals in the ring} $Z\left[\xi\right]$.\\
\bigskip\\
  \textbf{Theorem 1.2.}(The Reciprocity Law)( $\left[3\right] $). \textit{Let} $l$ \textit{be a prime odd natural number and} $\xi$ \textit{be a primitive root of l-order of unity. Let} $\alpha$ \textit{be a integer number}, $\alpha$ \textit{is not divisible with} $l$, $a$ \textit{be a semiprimary element from}$Z\left[\xi\right]$,$\alpha$ \textit{relatively prime with} $l$. 
  
    \textit{Then}:\\ 
   $$\left\{\frac{a}{\alpha}\right\}=\left\{\frac{\alpha}{a}\right\}.$$
  \bigskip\\
  \textbf{Theorem 1.3.}( $\left[ 3\right]$). \textit{Let} $\xi$ \textit{be a primitive root of l-order, of unity, where l is a prime natural number and let}
  $\mu$$\in$$Z\left[\xi\right]$.
  \textit{Let the Kummer field} $Q\left(\sqrt[l]{\mu};\xi\right).$ \textit{Then a prime ideal} $P$ \textit{in the ring} $Z\left[\xi\right]$, \textit{is in one of the cases}: \\
  (i)\textit{If} $\left\{\frac{\mu}{P}\right\}=0,$ \textit{then} $P$ \textit{is in the ring of integers A in the Kummer field} $Q\left(\sqrt[l]{\mu};\xi\right)$
\textit{equal with the l-power of a prime ideal;} \\
(ii) \textit{If} $\left\{\frac{\mu}{P}\right\}=1,$ \textit{then} $P$ \textit{ decomposes in l different prime ideals in the ring A;}\\
(iii) \textit{If} $\left\{\frac{\mu}{P}\right\}$ \textit{is equal with a root of order l of unity, different from} $1,$   \textit{then} $P$ \textit{is a prime ideal in the ring A.} 
\bigskip\\  
  \textbf{Proposition 1.4.}( $\left[10\right]$). \textit{Let A be the ring of integers of the Kummer field} $ \textbf{Q}\left(\sqrt[l]{p};\xi\right)$ \textit{where} $p$ \textit{is a prime natural number and} $\xi$ \textit{is a primitive root of order l of unity}. \textit{Let} $G$ \textit{be the Galois group of the Kummer field} $ \textbf{Q}\left(\sqrt[l]{p};\xi\right)$ \textit{over} \textbf{Q}. \textit{Then for any} $\sigma$$\in$$G$ \textit{and for any} $P$$\in$$Spec\left(A\right)$ \textit{we have} $\sigma\left(P\right)$$\in$$Spec\left(A\right)$.
\bigskip\\
\bigskip\\
 \section{Results}

\bigskip
  First, we state and prove two propositions that are necessary for solving  the equation\\ 
  $$x^{4}-q^{4}=py^{5}\ \ \     (1)$$\\
  in the conditions (2):\\
  (i) $p$ \textit{and} $q$ \textit{are different prime natural numbers;}\\
  (ii) $y$ \textit{is not divisible with} $p$;\\
  (iii)$\overline{p}$ \textit{is a generator of the group} $(U\left(Z_{q^{4}}\right),\cdot )$;\\
  (iv)$p\equiv 3$ ( \textit{mod} $20$ ), $q\equiv 4$ ( \textit{mod} $5$ );\\
  (v) gcd$\left(x,y\right)=1$;\\ 
  (vi) $2$ \textit{is a} $5$ \textit{-power residue mod} $q$.\ \ \ 
\bigskip\\  
\textbf{Lema 2.1. }\textit{Let} $p$ \textit{and} $q$ \textit{be prime integers satisfying the conditions} (2) \textit{and take} $\xi$ \textit{as a primitive root of order} $5$ \textit{of the unity}. \textit{If} $\textbf{Q}\left(\xi;\sqrt[5]{p}\right)$ \textit{is the Kummer field with the ring of integers A}, $y_{1}$ \textit{and} $y_{2}$ \textit{are integer numbers such that gcd}$\left(y_{1},y_{2}\right)=1$, $y^{5}_{2}-py^{5}_{1}=2q^{2}$, \textit{then, taking} $m$,$n$$\in$$\left\{0,1,...4\right\}$, $m$$\neq$$n$, \\
$$\left(y_{2}-\xi^{m}\sqrt[5]{p}y_{1}\right)A \  \textit{and}\ \left(y_{2}-\xi^{n}\sqrt[5]{p}y_{1}\right)A$$\\
 \textit{are comaximal ideals of A}.
\bigskip\\
\textbf{Lema 2.2. }\textit{Let us consider} $p$ \textit{and} $q$ \textit{satisfying the conditions} (2) \textit{and take} $\xi$ \textit{as a primitive root of order} $5$ \textit{of the unity}. \textit{If} $\textbf{Q}\left(\xi;\sqrt[5]{8p}\right)$ \textit{is the Kummer field with the ring of integers A}, $y_{1}$ \textit{and} $y_{2}$ \textit{are integers numbers, gcd}$\left(y_{1},y_{2}\right)=1$, $y^{5}_{2}-8py^{5}_{1}=q^{2}$, \textit{then, taking,} $m$,$n$$\in$$\left\{0,1,...,4\right\}$, $m$$\neq$$n$, \\
$$\left(y_{2}-\xi^{m}\sqrt[5]{8p}y_{1}\right)A \  \textit{and}\\ \left(y_{2}-\xi^{n}\sqrt[5]{8p}y_{1}\right)A$$\\
 \textit{are comaximal ideals of A}.
  \bigskip\\
  \textbf{Proof.} The proof of the Lema 2.1. and the proof of the the Lema 2.2. are similar with the proof of the Lema 2.1. from $\left[9\right]$.
  \bigskip\\
  Now we try to solve the equation $x^{4}-q^{4}=py^{5}$. 
  \bigskip\\
  \textbf{Theorem 2.3.} \textit{The equation} $x^{4}-q^{4}=py^{5}$ \textit{does not have nontrivial integer solutions in the conditions (2)}.
  \bigskip\\
  \textbf{Proof.} We suppose that the equation (1) has nontrivial integer solutions $\left(x,y\right)$$\in$$\textbf{Z}^{2}$ satisfying the conditions (2). We consider two cases: either \textit{x} is odd or \textit{x} is even.\\
  \textbf{Case I}: \textbf{x is an odd number}\\
  Knowing that \textit{q} is a prime natural number, $q\geq3$, we get $x^{2},q^{2}\equiv1$ (mod $4$) and therefore $x^{2}-q^{2}\equiv0$ (mod $4$), $x^{2}+q^{2}\equiv2$ (mod $4$).\\
 We denote $d=gcd\left(x^{2}-q^{2},x^{2}+q^{2}\right)$. Then $d/2x^{2}$ and  $d/2q^{2}$. But $gcd\left(x,y\right)=1$ implies $x$ is not divisible with $q$. Therefore $d=2$. We get either that\\
 \begin{center}
$x^{2}-q^{2}=16py^{5}_{1}$, $x^{2}+q^{2}=2y^{5}_{2},$
\end{center}
where $y_{1},y_{2}$$\in$\textbf{Z}, $2y_{1}y_{2}=y$, $y_{2}$ is an odd number, $gcd\left(y_{1},y_{2}\right)=1$ or that \\
\begin{center}
$x^{2}-q^{2}=16y^{5}_{1}$, $x^{2}+q^{2}=2py^{5}_{2},$
\end{center}
where $y_{1},y_{2}$$\in$\textbf{Z}, $2y_{1}y_{2}=y$, $y_{2}$ is an odd number, $g.c.d.\left(y_{1},y_{2}\right)=1.$\\
In the last case, we obtain that $p/\left(x^{2}+q^{2}\right)$, in contradiction with the fact that $p\equiv3$ (mod$4$).
It remains to study the case 
\begin{center}
$x^{2}-q^{2}=16py^{5}_{1}$, $x^{2}+q^{2}=2y^{5}_{2}$.
\end{center}
By substracting the two equations, we obtain $q^{2}=y^{5}_{2}-8py^{5}_{1}$.\\
Let A be the ring of integers of the Kummer field $\textbf{Q}\left(\xi;\sqrt[5]{8p}\right)$, where $\xi$ is a primitive root of order $5$ of unity. In A, the last equality becomes:\\ $$q^{2}=\left(y_{2}-y_{1}\sqrt[5]{8p}\right)\left(y_{2}-y_{1}\xi\sqrt[5]{8p}\right)...\left(y_{2}-y_{1}\xi^{4}\sqrt[5]{8p}\right).\ \ \ (3)$$
$p$$\equiv$$3$ (mod $5$). This implies (using the Proposition 1.1) that $p$ is a prime element in the ring $Z\left[\xi\right]$.\\
We prove that $\left\{\frac{q}{\left(p\right)}\right\}=1$
The conjugates of $\left\{\frac{q}{\left(p\right)}\right\}$ are $\left\{\frac{q}{\left(p\right)}\right\}^{2}$, $\left\{\frac{q}{\left(p\right)}\right\}^{3}$, $\left\{\frac{q}{\left(p\right)}\right\}^{4}$. But $\overline{\left\{\frac{q}{\left(p\right)}\right\}}=\left\{\frac{q}{\left(p\right)}\right\}$, therefore $\left\{\frac{q}{\left(p\right)}\right\}^{2}=\left\{\frac{q}{\left(p\right)}\right\}$. Knowing that $\left\{\frac{q}{\left(p\right)}\right\}\neq0$, we obtain that $\left\{\frac{q}{\left(p\right)}\right\}=1$.
   But $q\equiv4(mod 5)$ and, using the Proposition 1.1., we obtain that the ideal $qZ\left[\xi\right]=Q_{1}Q_{2}...Q_{r}$, $i=\overline{1,r}$, $r=\frac{\varphi\left(5\right)}{f}$ where $f=ord_{Z^{*}_{5}}\left(\overline{q}\right)=2$. Therefore $r=2$. Using the fact that $Z\left[\xi\right]$ is a principal ring, we obtain $qZ\left[\xi\right]=Q_{1}Q_{2}$, where $Q_{i}=\pi_{i}Z\left[\xi\right]$ and 
   $\pi_{i}$ are prime elements in the ring $Z\left[\xi\right]$, $i=1;2$.\\
   Now we try to decompose the ideal $(q)$ in the ring $A.$ For this we calculate:\\
   $\left\{\frac{8p}{\left(\pi_{i}\right)}\right\}=\left\{\frac{2}{\left(\pi_{i}\right)}\right\}^{3}\left\{\frac{p}{\left(\pi_{i}\right)}\right\}$ $i=1,2.$\\
 $2$ is a $5$-power residue mod $q$, implies there is $\alpha$$\in$$Z\left[\xi\right]$ such that $\alpha^{5}$$\equiv$$2$ (mod $q$), therefore $\alpha^{5}$$\equiv$$2$ (mod $\pi_{i}$), $i=1;2$. This implies that $\left\{\frac{2}{\left(\pi_{i}\right)}\right\}=1$, $i=1,2.$\\
 We obtain that $\left\{\frac{8p}{\left(\pi_{i}\right)}\right\}=\left\{\frac{p}{\left(\pi_{i}\right)}\right\}$, $i=1;2$.
   We have $1=\left\{\frac{q}{\left(p\right)}\right\}=\left\{\frac{\pi_{1}}{\left(p\right)}\right\}\cdot\left\{\frac{\pi_{2}}{\left(p\right)}\right\}=\left\{\frac{p}{\left(\pi_{1}\right)}\cdot\right\}\left\{\frac{p}{\left(\pi_{2}\right)}\right\}$ (according to the Theorem 1.7.)\\
   $\left\{\frac{p}{\left(\pi_{i}\right)}\right\}=\xi^{c_{i}}$$\equiv$$p^{\frac{N(\pi_{i})-1}{5}} (mod \pi_{i})$, $i=\overline{1,2}$.
   Using the fact that $q=\pi_{1}\pi_{2}$, where $\pi_{1},\pi_{2}$ are prime elements in the ring $Z\left[\xi\right]$ and that $N(q)=q^{4}$, we obtain that 
  $N(\pi_{1})=N(\pi_{2})=q^{2}$.\\
  If we suppose that $\left\{\frac{p}{\left(\pi_{1}\right)}\right\}=\left\{\frac{p}{\left(\pi_{2}\right)}\right\}=1$, it results that, for $i=1;2$, we have $p^{\frac{N(\pi_{i})-1}{5}}$$\equiv$ $1(mod \pi_{i})$, therefore $p^{\frac{q^{2}-1}{5}}$$\equiv$ $1(mod q)$. Since $\overline{p}$ is a generator of the group $\left(U\left(\textbf{Z}_{q^{4}}\right),\cdot\right)$ it results that $\varphi(q^{4})/\frac{q^{2}-1}{5}$. This implies that there exist $k$$\in$$\textbf{N}^{*}$ such that $\frac{q^{2}-1}{5}=k(q^{4}-q^{3})$. The last equality is equivalent with $q+1=5q^{3}k$. I obtained a contradiction with the fact $k,q$$\in$$\textbf{N}^{*}$, $q\geq3$.\\
  From the previously proved and from the fact that $1=\left\{\frac{p}{\left(q\right)}\right\}=\left\{\frac{p}{\left(\pi_{1}\right)}\right\}\cdot\left\{\frac{p}{\left(\pi_{2}\right)}\right\}$, we obtain that $\left\{\frac{p}{\left(\pi_{i}\right)}\right\}=\xi^{c_{i}}\neq1$, for $i=\overline{1,2}$. Using the Theorem 1.8., we obtain that $\pi_{1}A$, $\pi_{2}A$ are prime ideals in the ring A.\\
  Passing to ideals in the relation (3), we get:\\
   $$\left(y_{2}-y_{1}\sqrt[5]{8p}\right)A\left(y_{2}-y_{1}\xi\sqrt[5]{8p}\right)A...\left(y_{2}-y_{1}\xi^{4}\sqrt[5]{8p}\right)A=\left(\pi_{1}A\right)^{2}\left(\pi_{2}A\right)^{2}.\ \ \ (4)$$\\
   According to the Lema 2.2,the last equality is impossible.\\
   \textbf{The case II}: \textbf{x is an even number.}\\   
In this case, $x^{2}-q^{2}$ and $x^{2}+q^{2}$ are odd numbers.\\
We prove that $gcd\left(x^{2}-q^{2},x^{2}+q^{2}\right)=1$. We suppose that there exists an odd prime natural number $d$ such that $d/\left(x^{2}-q^{2}\right)$ and $d/\left(x^{2}+q^{2}\right)$. Hence $d/x$ and $d/q$. Using the hypothesis, we obtain that $d/y$, in contradiction with the fact $\left(x,y\right)=1$.  Therefore $gcd\left(x^{2}-q^{2},x^{2}+q^{2}\right)=1$.
Then (1) becomes: \\
$x^{2}-q^{2}=py^{5}_{1}$, $x^{2}+q^{2}=y^{5}_{2}$,
with $y_{1},y_{2}$$\in$\textbf{Z}, $y_{1}y_{2}=y$, $gcd\left(y_{1},y_{2}\right)=1$ or:\\
$x^{2}-q^{2}=y^{5}_{1}$, $x^{2}+q^{2}=py^{5}_{2}$,
with $y_{1},y_{2}$$\in$\textbf{Z}, $y_{1}y_{2}=y$, $gcd\left(y_{1},y_{2}\right)=1.$\\
In the last case, we obtain that $p/\left(x^{2}+q^{2}\right)$, in contradiction with the fact that $p\equiv3$ (mod$4$).
It remains to study the case 
\begin{center}
$x^{2}-q^{2}=py^{5}_{1}$, $x^{2}+q^{2}=y^{5}_{2}.$
\end{center}
Substracting the two equations, we get $2q^{2}=y^{5}_{2}-py^{5}_{1}$.\\
Let A be the ring of integers of the Kummer field $\textbf{Q}\left(\xi;\sqrt[5]{p}\right)$, where $\xi$ is a primitive root of order $5$ of the unity. In A, the last equality becomes:
$$\left(y_{2}-y_{1}\sqrt[5]{p}\right)\left(y_{2}-y_{1}\xi\sqrt[5]{p}\right)...\left(y_{2}-y_{1}\xi^{4}\sqrt[5]{p}\right)=2q^{2}.\ \ \ (5)$$ 
Similar with the case when $x$ is an odd number, we obtain $qA=\pi_{1}A\cdot\pi_{2}A$, where $\pi_{1}$,$\pi_{2}$ are irreducible elements in the rings $\textbf{Z}\left[\xi\right]$.\\
$2$ is a prime element in the ring $\textbf{Z}\left[\xi\right]$(according to Proposition 1.1.),
$p$ is a prime natural number, $p$$\equiv$$3(mod 20)$, implies $p$$\equiv$$1^{5}(mod 2)$, therefore $\left\{\frac{p}{\left(2\right)}\right\}=1$. Using the Theorem 1.8. we obtain that $2A=P_{1}P_{2}...P_{5}$, where $P_{1},P_{2},...,P_{5}$ are prime ideals in the ring $A$. \\
We consider the corresponding ideals in the relation (5) and we obtain:\\ 
$$\left(y_{2}-y_{1}\sqrt[5]{p}\right)A\left(y_{2}-y_{1}\xi\sqrt[5]{p}\right)A...\left(y_{2}-y_{1}\xi^{4}\sqrt[5]{p}\right)A=P_{1}P_{2}...P_{5}(\pi_{1}A)^{2}(\pi_{2}A)^{2}.\ \ \ (6)$$\\
Let G be the Galois group of the Kummer field $\textbf{Q}(\xi,\sqrt[5]{p})$ over \textbf{Q}. There is $v$$\in$G such that $v\left(\xi\right)=\xi$, $v\left(\sqrt[5]{p}\right)=\xi\sqrt[5]{p}$.
\smallskip\\
\textbf{Case (i):} If there exists $k$$\in$$\left\{1,2,...,5\right\}$ such that $\left(y_{2}-y_{1}\sqrt[5]{p}\right)A=P_{k}$$\in$Spec(A), we use Proposition 1.9., and we obtain that $v\left(\left(y_{2}-y_{1}\sqrt[5]{p}\right)A\right)=\left(y_{2}-y_{1}\xi\sqrt[5]{p}\right)A$$\in$Spec(A) and $v^{2}\left(\left(y_{2}-y_{1}\sqrt[5]{p}\right)A\right)=\left(y_{2}-y_{1}\xi^{2}\sqrt[5]{p}\right)A$$\in$Spec(A),..., $v^{4}\left(\left(y_{2}-y_{1}\sqrt[5]{p}\right)A\right)=\left(y_{2}-y_{1}\xi^{4}\sqrt[5]{p}\right)A$$\in$Spec(A),
therefore the equality (6) is impossible.
\smallskip\\
\textbf{Case (ii):} If there are $k$ and $h$ in $\left\{1,2,3,4,5\right\}$, $k\neq h$ such that
$\left(y_{2}-y_{1}\sqrt[5]{p}\right)A=P_{k}P_{h},$  $P_{k}$,$P_{h}$$\in$Spec(A), then we use Proposition 1.9. and Lema 2.1. and we obtain a contradiction with the relation (6)
\smallskip\\
\textbf{Case (iii):} If $\left(y_{2}-y_{1}\sqrt[5]{p}\right)A=\left(\pi_{1}A\right)^{2}$, then (in according to the Proposition 1.9.) $\left(y_{2}-y_{1}\xi\sqrt[5]{p}\right)A=\left(\pi_{2}A\right)^{2}$,..., $\left(y_{2}-y_{1}\xi^{4}\sqrt[5]{p}\right)A=P^{2}$, $P$$\in$Spec(A), in contradiction with (6).\\
  From cases (i), (ii), (iii), it results that the equality (6) is impossible.
  We get that the equation (1) does not have nontrivial integer solutions satisfying the conditions (2).
  \bigskip\\
  \textbf{References}
  \bigskip\\
  $\left[1\right]$ T. Albu, I.D. Ion, \textit{Chapters of the algebraic theory of numbers} (in Romanian), Ed. Academiei, Bucuresti, 1984.\\ 
  $\left[2\right]$ V. Alexandru, N.M. Gosoniu, \textit{Elements of Numbers Theory} (in Romanian), Ed. Universitatii, Bucuresti, 1999.\\ 
  $\left[3\right]$ D.Hilbert, \textit{The theory of algebraic number fields}, Ed. Corint, Bucuresti, 1998.\\
  $\left[4\right]$ K.Ireland, M.Rosen \textit{A Classical Introduction to Modern Number Theory }, Springer-Verlag, 1992.\\
  $\left[5\right]$ D. Savin, \textit{On some Diophantine Equations (I)}, An. St. Univ. "Ovidius", Constanta, Ser. Mat., \textbf{10} (2002), fasc.1, p.121-134.\\
  $\left[6\right]$ D. Savin, \textit{On some Diophantine Equations (II)}, An. St. Univ. "Ovidius", Constanta, Ser. Mat., \textbf{10} (2002), fasc.2, p.79-86.\\
  $\left[7\right]$ D. Savin, \textit{Systems of Diophantine Equations without solutions}, Proceedings of the $10^{th}$ International Symposium of Mathematics and its Applications, "`Politehnica"' University of Timisoara, November 6-9, 2003, p. 310 - 317.\\
  $\left[8\right]$ D. Savin, \textit{On some Diophantine Equations (III)}, An. St. Univ. "Ovidius", Constanta, Ser. Mat., \textbf{12}(2004), fasc.1, p. 73 - 80.\\
  $\left[9\right]$ D. Savin, \textit{On the Diophantine Equation} $x^{4}-q^{4}=py^{3}$, An. St. Univ. "Ovidius", Constanta, Ser. Mat., \textbf{12}(2004), fasc.1, p. 81 - 90.\\
  $\left[10\right]$ M. Stefanescu, \textit{Galois Theory} (in Romanian), Ed. Ex Ponto, Constanta, 2002.\\ 
  $\left[11\right]$ C. Vraciu, M. Vraciu, \textit{Basics of Arithmetic} ( in Romanian), Ed. All, Bucuresti, 1998.\\ 
  \bigskip\\

  {Department of Mathematics\\
"Ovidius" University of Constanta}\\
Bd. Mamaia 124, Constanta\\
Romania\\
{e-mail: Savin.Diana@univ-ovidius.ro\\
dianet72@yahoo.com} 
 \end{document}